\documentclass[a4paper]{article}

\usepackage{amsmath, amssymb, amsfonts, amscd, amsthm}

\newcommand{\Z}{\mathbb{Z}}
\newcommand{\C}{\mathbb{C}}
\renewcommand{\c}{}
\newcommand{\R}{\mathbb{R}}
\newcommand{\Q}{\mathbb{Q}}
\renewcommand{\H}{\mathbb{H}}

\newcommand{\F}{\mathbb{F}}

\newcommand{\<}{\langle}
\renewcommand{\>}{\rangle}

\newcommand{\af}{\mathbf{A}}
\newcommand{\ga}{\mathbf{G}_a}
\newcommand{\gm}{\mathbf{G}_m}
\newcommand{\gc}{\mathbf{G}_2}
\newcommand{\SO}{\mathbf{SO}}
\renewcommand{\O}{\mathbf{O}}
\newcommand{\SL}{\mathbf{SL}}
\newcommand{\GL}{\mathbf{GL}}
\newcommand{\Sp}{\mathbf{Sp}}
\newcommand{\Spin}{\mathbf{Spin}}
\newcommand{\T}{T}

\newcommand{\cht}[1]{CH^*#1 \otimes_\Z \F_2}
\newcommand{\chp}[1]{CH^*#1 \otimes_\Z \F_p}
\newcommand{\ich}{\mathfrak{ch}}

\newcommand{\VR}{V_\mathbb{R}}
\newcommand{\VC}{V_\mathbb{C}}
\newcommand{\DR}{\Delta_\mathbb{R}}
\newcommand{\DC}{\Delta_\mathbb{C}}
\newcommand{\D}{\Delta}
\newcommand{\z}{\textrm{\bf z}}
\newcommand{\x}{\textrm{\bf x}}
\newcommand{\y}{\textrm{\bf y}}
\renewcommand{\zeta}{\xi}

\newcommand{\td}[1]{\widetilde{#1}}

\newtheoremstyle{pedro}{}{}{\itshape}{}{\sc}{~--}{ }{\thmname{#1}\thmnumber{ #2}\thmnote{ (#3)}}

\theoremstyle{pedro}
\newtheorem{lem}{Lemma}[section]

\newtheorem{thm}[lem]{Theorem}
\newtheorem*{thm2}{Theorem}

\newtheorem{prop}[lem]{Proposition}

\newtheorem{coro}[lem]{Corollary}

\theoremstyle{remark}

\newtheorem{rmk}[lem]{Remark}
\newtheorem*{rmk2}{Remark}

\theoremstyle{definition}

\usepackage{titlesec}

\titleformat{\section}{\bf\center}{\S \arabic{section}.}{1em}{}
\titleformat{\subsection}{\bf}{}{0pt}{}   

\title{The Chow rings of $G_2$ and $Spin(7)$}

\author{Pierre Guillot}

\begin{document}

\maketitle
\begin{abstract} 
We compute $CH^*B\gc$ and $CH^*(B\Spin_7\c)_{(2)}$, using in
particular the stratification method introduced by Vezzosi. We also
give some information on the Chow rings of the finite groups $G_2(q)$.

These results together with work by Schuster and Yagita imply the
existence of a counter-example to a conjecture of Totaro's.

\smallskip
{\bf Status:} this paper has been published in J. reine
angewand. Math. vol 604, 2007. This version includes an appendix which
will appear separately, also in J. reine angewand. Math. It is in this
appendix that Totaro's conjecture is mentioned. It does not appear in
the summary of results in the Introduction, for we want this document
to remain as close as possible to the published version.

\end{abstract}

\section{Introduction.}

Since Totaro proved (\cite{totaro}) that the classifying space $BG$ of an algebraic group $G$ could be approximated by algebraic varieties, a large number of papers have appeared on the subject of computing the Chow ring $CH^*BG$. Totaro computed this himself for $\GL_n$, $\SL_n$, $\Sp_n$, $\O_n$ and $\SO_{2n+1}$; see the work of Pandharipande on $\SO_4$ \cite{panda}; Rebecca Field computed $CH^*B\SO_{2n}$ \cite{reb}; Vezzosi worked on $\mathbf{PGL}_3$ \cite{vezzosi} while Vistoli has subsequently obtained results on $\mathbf{PGL}_p$ \cite{vistoli2}; a number of these results are given alternative proofs by Molina and Vistoli in \cite{vistoli1}; see also the work of Yagita \cite{yag1}, \cite{yag2}, \cite{yag3}. The author has computed the mod $p$ Chow rings of $BG$ for certain finite groups $G$ (symmetric groups and Chevalley groups in particular), see \cite{pedro1} and \cite{pedro2}. We refer to all these papers for more background and motivation.

The first object of the present article is to compute the Chow ring of the exceptional group $\gc$. This is the first complete calculation with exceptional groups. We obtain the following (theorem \ref{thm:main} in the text).

\begin{thm2} Let $c_i$ denote the $i$-th Chern class of $V$, the $7$-dimensional irreducible representation of the exceptional group $\gc$. Then one has
$$CH^*B\gc=\Z[c_2,c_4,c_6,c_7]/(c^2_2=4c_4, c_2c_7=0,2c_7=0).$$
\end{thm2}

In a series of paper (\cite{yag1}, \cite{yag2}, \cite{yag3}), Yagita has obtained the localised Chow ring $CH^*B\gc \otimes_\Z \Z_{(2)}$. We point out that our proof is more elementary. Among other things, Yagita relies on the proof of the Bloch-Kato conjecture for $p=2$ by Voevodsky, on motivic cobordism, and on the Atiyah-Hirzebruch motivic spectral sequence. By contrast, we use a theorem of Totaro's (\ref{thm:totaro} below) which relates an algebraic group $G$ to the normalizer of a maximal torus, and our proof is mainly a study of two subgroups of maximal rank in $\gc$, together with a few elementary remarks on representation theory.

Our second result concerns the finite groups $G_2(q)$. In previous work cited above, we have obtained some information on the mod $p$ Chow rings of these, for $p\ge 3$. Here we complete the picture by proving the following (see \ref{thm:g2q} and corollary).

\begin{thm2} Suppose $q \equiv 1$ mod $4$. Then the Chow ring $\cht{BG_2(q)}$ contains a polynomial subring $\F_2[d_4,d_6,d_7]$, which injects into the mod $2$ cohomology ring. Unlike the situation when $p$ is odd, the restriction 
$$\cht{BG_2(q)} \to \cht{BT(q)}$$ is not injective, where $T$ is a maximal torus.
\end{thm2}

This result does not seriously depend on the previous one, although it comes as a natural complement.

Our final and main result concerns the Chow ring of the spinor group $\Spin_7\c$, of which $\gc$ is a subgroup. This is the first "genuine" spin group, since for each $n\le 6$, $\Spin_n\c$ is isomorphic to a group of type $A$ or $C$ (or a product of such). After proving that $CH^*B\Spin_7\c$ has only $2$-torsion (\ref{lem:torsionspin}), we obtain a nearly complete calculation after localizing at $2$. The theorem below is a summary of corollary \ref{coro:generators}, corollary \ref{coro:addspin}, and proposition \ref{prop:relations}.

\begin{thm2} Let $c_i$ denote the $i$-th Chern class of the $7$-dimensional representation of $\Spin_7\c$, and let $c'_i$ denote the $i$-th Chern class of the spin representation. Then there is an additive isomorphism $$CH^*(B\Spin_7\c)_{(2)}= \Z_{(2)}[c_4,c_6,c'_8]\otimes \Big( \Z_{(2)}\<1,c'_2,c'_4,c'_6\> \oplus \F_2\< \zeta_3 \> \oplus \F_2[c_7]\<c_7\> \Big)$$ where $\zeta_3$ is a class of (cohomological) degree $6$, which cannot be expressed in terms of transfers of Chern classes. The products in $CH^*(B\Spin_7\c)_{(2)}$ are determined by the following relations.
$$\begin{array}{rlrl}
(1) & \zeta_3^2=0 & (8) & c'_2c'_6 - c'_2c_6=\frac{2}{3}c_4(c'_4 - c_4) + 16c'_8  
\\ (2) & \zeta_3c_7=0 &  (9) & c'_2c_7=\delta_1c_6\zeta_3
\\ (3) & \zeta_3c'_4=\zeta_3c_4 & (10) &  c'_4(c'_4 - c_4)=c_4(c'_4 - c_4) + 36 c'_8  
\\ (4) & \zeta_3c'_6=\zeta_3c_6 & (11) & c'_4(c'_6 - c_6)=c_4(c'_6 - c_6) +6c'_2c'_8  
\\ (5) & \zeta_3c'_2=0 & (12) & c'_4c_7 - c_4c_7=\delta_2c_8\zeta_3 
\\ (6) & (c'_2)^2 - 4c_4=\frac{8}{3}(c_4'-c_4) & (13) & c'_6(c'_6 - c_6)=c_6(c'_6 - c_6) + c'_8(\frac{8}{3}c'_4 + \frac{4}{3}c_4)
\\ (7) & c'_2c'_4 - c'_2c_4 = 6(c'_6 - c_6) & (14) & c'_6c_7 - c_6c_7=0
\end{array}$$ Here $\delta_i=0$ or $1$ ($i=1,2)$

In particular, the restriction map to $CH^*(B\gc)_{(2)}$ is surjective with kernel generated by $\zeta_3$, $c'_4 - c_4$, $c'_6 - c_6$ and $c'_8$.

Finally, there is an isomorphism
$$CH^*(B\Spin_7\c)_{(2)}=BP^*(B\Spin_7\c)\otimes_{BP^*} \Z_{(2)}$$ where $BP$ is the Brown-Peterson spectrum at the prime $2$.
\end{thm2}

This result completes thus the computations of Kono-Yagita on the additive $BP$-cohomology of $\Spin_7\c$ in \cite{konoyagita}.

It may come as a surprise to someone acquainted with the known calculations in the field that it should take $16$ relations (counting $2\zeta_3=0$ and $2c_7=0$) to describe the Chow ring of $\Spin_7\c$. Let us consider mod $2$ Chow rings to simplify the discussion. Then $\GL_n$, $\SL_n$, $\O_n$, $\SO_{2n+1}$, and $\Sp_n$ have a polynomial Chow ring mod $2$, and the situation for $\SO_{2n}$ is not much worse (see theorem \ref{thm:reb} below). On the other hand, $\Spin_7$ offers a much more involved answer.

Likewise, Chow rings tend to be "simpler" algebraically than the corresponding cohomology rings (this is particularly true for finite groups (see \cite{pedro2}) but remains valid for connected groups, $\O_n$ being a good example). Here, the mod $2$ cohomology of $\Spin_7$ is polynomial, so we have a strong exception to the rule.

\subsection{Organization of the paper.}

This work is divided into two parts. The first is dedicated to $G_2$: we give its definition, study two subgroups of maximal rank, and investigate the Chern classes of representations in sections \ref{section:intro}, \ref{section:maximalrank} and \ref{section:chern} respectively. The computation of $CH^*B\gc$ is completed at this point. Then in \S \ref{section:g2q} we make a few comments on the finite groups $G_2(q)$. This section is not needed to understand the sequel, unlike the beginning.

Part II is devoted to $\Spin_7\c$. We study its action under the spin representation in \S \ref{section:prelim}, in order to implement the "stratification method" in \S \ref{section:stratspin}. We make use of Brown-Peterson cohomology in \S \ref{section:otherinput}, and we put all the ingredients together in \S \ref{section:end}.

\subsection{Acknowledgements.}

It is a pleasure to thank Angelo Vistoli for his early interest in this work, and for many helpful discussions. Nobuaki Yagita has helped me clarify my understanding of his computations, and I would like to express my gratitude to him. Mike Shuter has run many computer calculations using Magma which have been useful in checking my results -- I am sorry that they have not made their way into the finalized version of this paper, but I am very grateful for them.

Finally, it is not the first time, and certainly not the last, that I find myself thanking Burt Totaro warmly for a wealth of comments and suggestions.

\subsection{Notations.}

We shall explain most of the notations as we introduce them. However, we point out here that all references to algebraic geometry are to be understood over the complex numbers.

Note also that we will write $SO(n)$, $SU(n)$ or $Spin(n)$, say, for compact Lie groups, while algebraic groups will be named $\SO_n\c$, $\GL_n\c$ or $\Spin_n\c$. As the Chow ring $CH^*BG$ depends on the structure of $G$ as an algebraic group, and not just on the underlying Lie group of complex points $G(\C)$, it may be important to distinguish, for instance, between $\SO_n\c$ and $\SO_n(\C)$. However we have not pushed this to the extent of writing $\boldsymbol{\mu}_2$ for the cyclic group $\Z/2$.

\pagebreak
\bigskip
\begin{center}
{\sc\Large Part I: The Chow ring of $G_2$ and related finite groups}
\end{center}

\section{The group $G_2$.}\label{section:intro}

\subsection{Definitions.}(See \cite{adams}.) A {\em normed algebra} over the reals is a finite dimensional, real vector space $A$ with a product $\mu:A\otimes A \to A$ and a norm satisfying $\|\mu(x,y)\|=\|x\| \cdot \|y\|$. The dimension of $A$ can only be $1$, $2$, $4$ or $8$, and in each respective case, there is only one possibility for $A$ up to isomorphism: the field $\R$ of real numbers, the field $\C$ of complex numbers, the (non-commutative) field of quaternions $\H$, and the non-associative algebra of Cayley numbers or octonions $\mathbb{O}$.

The group $G_2$ is defined to be the group of automorphisms of $\mathbb{O}$ as an algebra (the uniqueness of $\mathbb{O}$ implies that an algebra automorphism automatically preserves norms). It is a compact, $1$-connected Lie group of dimension $14$ and rank $2$. It is one of the five so-called exceptional Lie groups: these are the simple, $1$-connected compact Lie groups which are not isomorphic to either $SU(n)$, $Spin(n)$ or $Sp(n)$.

Let $\Delta_\R$ denote the spin representation of $Spin(7)$: it has real dimension $8$ and we may see it as the vector space underlying $\mathbb{O}$. In this way, $G_2$ can be seen as the subgroup of $Spin(7)$ fixing the unit $1\in \mathbb{O}$, and the action of $G_2$ on $\mathbb{O}$ is via $Spin(7)$.

The algebraic group $\gc$ is the complexification of $G_2$ (see \cite{btd}, III, 8). One possible definition is thus as follows: consider the inclusions $G_2\subset SO(7) \subset \SO_7\c \subset \GL_7\c$, and define $\gc$ to be the Lie subgroup of $\GL_7\c$ whose Lie algebra is the complexification $\mathfrak{g}_2 \otimes_\R \C$, where $\mathfrak{g}_2$ is the (real) Lie algebra of $G_2$. However, one may also consider the complex algebra $\mathbb{O} \otimes_\R \C$ and define $\gc$ to be its group of automorphisms as an algebra. In any case, $G_2$ is a maximal compact subgroup of $\gc$, and we have homotopy equivalences $G_2\simeq \gc$ and $BG_2\simeq B\gc$.

The group $G_2$ comes equipped with a $7$-dimensional representation which we denote $\VR$, via $G_2 \to Spin(7) \to SO(7)$. In fact, the spin representation $\Delta_\R$, when restricted to $G_2$, splits as the direct sum of $\VR$ and a trivial $1$-dimensional representation. Thus we may also regard $\VR$ as the action of $G_2$ on "purely imaginary" Cayley numbers, by which we mean the orthogonal complement of $1$.
The complexification of $\VR$ is denoted $\VC$ or simply $V$. It extends to a representation of $\gc$.

\subsection{Known information on the cohomology and Chow ring of $G_2$.} Throughout this article, we let $w_i$, resp. $c_i$, denote the $i$-th Stiefel-Whitney class of $\VR$, resp. the $i$-th Chern class of $\VC$. Recall that $c_i=w_i^2$ in mod $2$ cohomology (\cite{ms}, problem 14-B).


It is then known that (\cite{borel1}):
$$H^*(B\gc,\F_2)=H^*(BG_2,\F_2)=\F_2[w_4,w_6,w_7].$$
The rational cohomology, on the other hand, is:
$$H^*(B\gc,\Q)=CH^*B\gc\otimes_\Z \Q=\Q[c_2,c_6]$$ where the first isomorphism follows from \cite{edidin} (it would hold with $\gc$ replaced by any connected algebraic group); moreover, this ring is also isomorphic to $(H^*(BT,\Q))^W$ where $T$ denotes a maximal torus and $W$ is the Weyl group.

As for the integral Chow ring, Totaro proved that $CH^*B\gc$ is generated by the classes $c_i$, see \cite{totaro}. One of the aims of this paper is to determine the relations between these. A priori, one only knows that $2c_i=0$ when $i$ is odd, since $V$ has a real structure (equivalently, $V$ is self-dual), and that $c_1=0$, for example because $\gc$ acts on $V$ via $\SL_7\c$.

\section{Subgroups of maximal rank}\label{section:maximalrank}

\subsection{Torus normalizers.} 
The following will be a crucial ingredient in our computations. A
sketch of proof may be found in \cite{vezzosi}.

\begin{thm}[Totaro]\label{thm:totaro} 
Let $G$ be a linear algebraic group, and let $N$ denote the normalizer
of a maximal torus. Then the restriction map
$$CH^*BG \to CH^*BN$$ is injective.
\end{thm}

Given the strength of this result, it is natural to look for subgroups
of $\gc$ containing a maximal torus.

\subsection{The subgroup $\SL_3\c\rtimes \Z/2$.} The first example is as follows.

\begin{lem} The group $\gc$ contains a semi-direct product $K=\SL_3\c\rtimes \Z/2$.
\end{lem}

\begin{proof} 
Recall that $G_2$ acts on the sphere $S^6$ with stabilizers
$SU(3)$. Pick an element $x$ on $S^6$ and $g\in G_2$ which sends $x$
to $-x$. Then $g^2$ belongs to the stabilizer $S$ of $x$, so that
together $S$ and $g$ generate an extension of $SL_3\c$ by
$\Z/2$. Moreover, the centre of $SL_3\c$ is $\Z/3$ and
$H^2(\Z/2,\Z/3)=0$, so the extension has to be split, and $G_2$
contains a semi-direct product $SU(3)\rtimes \Z/2$. The result follows
by taking complexifications.
\end{proof}

\begin{coro} If a class $x\in CH^*B\gc$ is torsion, then $2x=0$.
\end{coro}
\begin{proof} Clearly, $K$ possesses an element which normalizes the standard maximal torus in $\SL_3\c$, but does not belong to $\SL_3\c$. Since the Weyl group of $\gc$ has order $12$ and that of $\SL_3\c$ has order $6$, it follows that $K$ contains the full normalizer (in $\gc$) of a maximal torus $\T$. (Recall that $\gc$ is reductive, so $\T$ is its own centraliser.) By theorem \ref{thm:totaro} above, the map $CH^*B\gc\to CH^*BK$ is injective.

Now, the Chow ring of $\SL_3\c$ is torsion-free and this group has
index $2$ in $K$, so the result follows from the projection formula.
\end{proof}

\subsection{A subgroup isomorphic to $\SO_4\c$.} For a proof that $G_2$ also contains a copy of $SO(4)$, we refer to \cite{mimuratoda}, II.5.7; see also \cite{borel2} and \cite{borel1}. We shall say more in remark \ref{rmk:so4} below.

\begin{lem}\label{lem:injectso4} The restriction map
$$CH^*B\gc \to CH^*B\SO_4\c$$ is injective.
\end{lem}
\begin{proof} Let $\T$ be a maximal torus in $\SO_4\c$, let $N'=N_{\SO_4\c}(\T)$ and $N=N_{\gc}(\T)$. Then $N'$ has index three in $N$, so that any $x\in CH^*BN$ restricting to $0$ in $CH^*BN'$ satisfies $3x=0$, and by theorem \ref{thm:totaro}, the same is true for an $x$ taken in $CH^*B\gc$, restricting to $0$ in $CH^*BN'$. By the corollary above, we must also have $2x=0$, so that $x=0$, and the restriction map to $N'$ is injective. The result follows.
\end{proof}

\section{Chern classes}\label{section:chern}

\subsection{Representations of $\SO_{2m}\c$.} We recall briefly a few facts concerning the representation theory of $\SO_{2m}\c$. See \cite{adams2} or better, \cite{milnor}. Let $W$ denote the standard representation, and let $\lambda_m=\Lambda^m(W)$ denote its $m$-th exterior power. Then $\lambda_m$ splits as the sum $\lambda_m^+ \oplus \lambda_m^-$ of the eigenspaces of the star operator (the corresponding eigenvalues being $\pm1$ or $\pm i$ according as $m$ is even or odd).

To describe the character of $\lambda_m^\pm$, we introduce the following notation. Let $\T$ denote the standard maximal torus of $\SO_{2m}\c$, let $\alpha_i$ denote the $i$-th standard character of $T$, and let $\tau_i$ denote the $i$-th elementary symmetric function in the $m$ variables $\alpha_j + \alpha_j^{-1}$. The function $\tau_m$ splits as $\tau_m=\tau^+_m + \tau^{-1}_m$ where $$\tau_r^\pm=\sum_{\epsilon_1\ldots\epsilon_m=\pm 1} \alpha_1^{\epsilon_1}\ldots\alpha_m^{\epsilon_m}.$$ 
Then one has:
$$\chi(\lambda_m^\pm)=\tau_m^\pm + \tau_{m-2} + \ldots + \frac{1}{2} {2k \choose k} \tau_{m-2k} + \ldots .$$

The Chern classes of $\lambda_m^+$ will play an essential role in the case at hand, namely when $m=2$. Indeed, consider the $7$-dimensional representation $V$ of $\gc$. The weights for $V$ are classically given as $\{0, \pm x_1, \pm x_2, \pm x_3 \}$, where the Lie algebra tangent to a maximal torus in $G_2$ is seen as the subspace $x_1 + x_2 + x_3 = 0$ in a $3$-dimensional vector space having dual basis given by the $x_i$'s. In other words, the character of $V$ restricted to $T$ is
$$1 + \alpha_1 + \alpha_1^{-1} + \alpha_2 + \alpha_2^{-1} + \alpha_1\alpha_2 + \alpha_1^{-1}\alpha_2^{-1}=1 + \tau_1 + \tau_2^+.$$
Since the character of $W$ is $\tau_1$, this shows: 
\begin{lem}\label{lem:vasso4}As an $\SO_4\c$-module, $V=W\oplus \lambda_2^+$.
\end{lem}

\begin{rmk}\label{rmk:so4} The above splitting may in fact be realised over $\R$. The summand corresponding to $\lambda_2^+$ is $\H$, the algebra of quaternions, while $W$ corresponds to the orthogonal complement $\H^\perp$. Now, $SO(4)$ has maximal rank in $G_2$, and looking at the Weyl groups shows that there is no connected group $G$ with $SO(4)\subsetneq G \subsetneq G_2$. Thus we see {\em a posteriori} that this $SO(4)$ may be defined as the subgroup of $G_2$ preserving $\H$ (which can be seen directly to be connected).

This remark is due to Burt Totaro.
\end{rmk}

\subsection{The Chow ring of $B\SO_4\c$.} As above, $\T$ is a maximal torus in $\SO_{2m}\c$. Write
$$CH^*B\T=\Z[t_1,\ldots,t_m]=H^*(B\T,\Z).$$
(Note that this is the only occasion in this paper where subscripts do not refer to the degree; each $t_i$ has (cohomological) degree $2$.)

Let us now recall the computation of the Chow ring of $\SO_{2m}\c$.

\begin{thm}[Field \cite{reb}]\label{thm:reb} The Chow ring of $\SO_{2m}\c$ is
$$CH^*B\SO_{2m}\c=\Z[d_2,\ldots,d_{2m},y_m]/(y_m^2=(-1)^m2^{2m-2}d_{2m},2d_{odd}=0, y_md_{odd}=0)$$ where the $d_i$'s are the Chern classes of $W$ and $y_m$ is a class restricting to $\pm 2^{m-1}t_1\ldots t_m$ in $CH^*BT$.
\end{thm}

When $m=2$, all these classes can be expressed in terms of Chern classes. Indeed, the total Chern class for $W$ in $CH^*B\T$ is
$$\sum d_iX^i=(1+t_1X)(1-t_1X)(1+t_2X)(1-t_2X)=1 - (t_1^2 + t_2^2)X^2 + t_1^2t_2^2X^4$$ so that, after restricting to $T$, one has $d_2=-(t_1^2 + t_2^2)$. As for $\lambda_2^+$, the total Chern class is simply
$$(1+(t_1+t_2)X)(1-(t_1+t_2)X)=1 - (t_1+t_2)^2X^2$$ and thus the second Chern class $c_2(\lambda_2^+)=d_2 + y_2$ after restricting to $\T$ (where we have chosen the sign of $y_2$ appropriately).

However, the theorem above ensures that $CH^2B\SO_4\c$ is torsion-free, while the restriction map $CH^*BG \to CH^*B\T$ is always injective modulo torsion, for any linear algebraic group $G$ with maximal torus $\T$. Therefore the relation $c_2(\lambda_2^+)=d_2 + y_2$ holds in $CH^*B\SO_4\c$.

\subsection{The Chow ring of $\gc$.} We can now state and prove our main result.

\begin{thm}\label{thm:main} Let $c_i$ denote the $i$-th Chern class of $V$, the $7$-dimensional irreducible representation of the exceptional group $\gc$. Then one has
$$CH^*B\gc=\Z[c_2,c_4,c_6,c_7]/(c^2_2=4c_4, c_2c_7=0,2c_7=0).$$
\end{thm}

\begin{proof} We know that the Chow ring of $\gc$ is generated by the $c_i$'s. By lemma \ref{lem:injectso4}, we may prove that the relations above hold by checking them in the Chow ring of $\SO_4\c$.

Now, by lemma \ref{lem:vasso4}, we obtain the relations below, where $d_3'$ denotes the third Chern class of $\lambda_2^+$ in the Chow ring of $\SO_4\c$ (we have seen above that $d'_3$ is $0$ in $CH^*B\T$, but it need not be $0$).
$$c_1=0,\quad c_2=2d_2 + y_2, \quad c_3=d_3 + d_3',\quad c_4=d_2^2 + d_2y_2 + d_4,$$ $$c_5=d_2d_3+d_2d_3',\quad c_6=d_3d_3' + d_2d_4 + y_2d_4,\quad c_7=d_4d_3'.$$ 
Now, if $d_3'$ were $0$, we would have $c_7=0$, but this class maps to $w_7^2\ne 0$ in $H^*(B\gc,\F_2)$. So $d_3'\ne 0$, and therefore we must have $d_3'=d_3$. It follows immediately that $c_3=c_5=0$, and that all the relations stated in the theorem hold.

We prove now that these are the only relations. In other words, we consider an element $P$ in the polynomial ring $\Z[C_2,C_4,C_6,C_7]$ such that $P(c_2,c_4,c_6,c_7)=0$ and we show that $P$ belongs to the ideal $I$ generated by $C_2^2 - 4C_4$, $C_2C_7$ and $2C_7$.

Write $P=Q + R$ where the variable $C_2$ does not appear in $R$ (and appears in each monomial of $Q$). Since $c_2$ maps to $0$ in mod $2$ cohomology, the relation $P(c_2,c_4,c_6,c_7)=0$ becomes $R(c_4,c_6,c_7)=0$ in $H^*(B\gc,\F_2)=\F_2[w_4,w_6,w_7]$. As $c_i=w_i^2$, we conclude that the coefficients of $R$ are $0$ mod $2$. Now, since $2C_7\in I$, we conclude that, modulo the ideal $I$, we may assume that $R$ does not contain the variable $C_7$.

Similarly, since $C_2C_7\in I$, we may assume mod $I$ that $Q$ does not contain the variable $C_7$. So we only need to consider the case where $P$ does not contain the variable $C_7$ at all.

Write the long division $P=A(4C_4 - C_2^2) + B$, where $A$ and $B$ are polynomial in the variables $C_2$, $C_4$ and $C_6$ with rational coefficients such that the degree of $B$ in $C_4$ is $0$ (ie $B$ does not contain the variable $C_4$). Since $CH^*B\gc\otimes \Q=\Q[c_2,c_6]$ is a polynomial ring, it follows that $B=0$, so that $C_2^2 - 4C_4$ divides $P$ in $\Q[C_4,C_4,C_6]$. Since the coefficient of $C_2^2$ in $C_2^2 -4C_4$ is $1$, we may write another long division and conclude that $P$ is in fact divisible by $C_2^2 -4C_4$ in $\Z[C_4,C_4,C_6]$.
\end{proof}

\begin{rmk} Working modulo $2$, we have
$$\cht{B\gc}=\F_2[c_2,c_4,c_6,c_7]/(c^2=0,c_2c_7=0)$$ and it follows that the cycle map to mod $2$ cohomology is not injective: indeed its kernel is the ideal generated by $c_2$. Note that $c_2$ restrict to $0$ in the (integral) Chow ring of any elementary abelian $2$-subgroup.
\end{rmk}

\section{The groups $G_2(q)$.}\label{section:g2q}

This section is a bit of a digression. The reader who is mainly interested in the Chow ring of $\Spin_7\c$ may skip it altogether. Nevertheless, we feel it appropriate to insert these results here, as they give an idea of the behavior of the whole $G_2$ family, as far as the ring of algebraic cycles in concerned.

 The complex algebraic group $\gc$ can be defined over $\Z$, and thus for any ring $R$ there is a group $\gc(R)$. In this section we shall be interested in the case of $\gc(k)$ where $k$ is a finite field, which in group theory is traditionnally denoted by $G_2(q)$ where $q$ is the order of $k$. Here to avoid confusion, we caution that we see $G_2(q)$ as an algebraic group {\em over the complex numbers}, as we may do for any finite group.

A definition is $G_2(q)$ is given in \cite{dickson}. In what follows, we shall only need information on its $2$-local structure, ie its $2$-Sylow subgroups and its elementary abelian subgroups. 

We let $T$ denote a maximal torus (defined and split over $\Z$) of the group scheme $G_2$, we let $N$ be its normalizer, and finally we let $W$ denote the Weyl group. Thus for any field $k$, we have
$W=N(k)/T(k)$ and when $k$ is algebraically closed, $N(k)$ is indeed the normalizer in the group-theoretic sense of $T(k)$ (see \cite{pedro2}, 3.2).

We have investigated the Chow rings of Chevalley groups, of which $G_2(q)$ is an example, in \cite{pedro2}. There we establish the following.

\begin{prop} Let $p$ be a prime number. Suppose $k$ is a finite field of characteristic different from $p$, containing the $p$-th roots of unity. Then if $p>3$, we have $$\chp{BG_2(k)}=(\chp{BT(k)})^W=\F_p[s_2,s_6].$$ If $p=3$, we have  at least 
$$\ich^*BG_2(k)=(CH^*BT(k)\otimes_\Z \F_3)^W=\F_3[s_2,s_6]$$
where $\ich^*BG_2(k)=im(CH^*BG_2(k)\otimes_\Z \F_3 \to H^*(BG_2(k),\F_3)$.
\end{prop}

The situation when $p=2$ is quite different, and we turn to this now. We shall consider a field $k=\F_q$ with $q\equiv 1$ mod $4$.

\subsection{The C\'ardenas-Kuhn theorem.} Let $L\subset K \subset G$ be finite groups. They are said to form a {\em closed system} if every subgroup of $K$ which is conjugate to $L$ in $G$ is also conjugate to $L$ in $K$. We let $N_G(L)$ denote the centralizer of $L$ in $G$, and $W_G(L)=N_G(L)/L$; similarly for $W_K(L)$. In the following statement, $h^*$ denotes a "cohomological inflation functor" or "cohomological global Mackey functor": see \cite{webb1} and \cite{webb2}. Roughly, what is meant is any contravariant functor from finite groups to, say, graded $\F_2$-vector spaces, endowed with transfers from subgroups satisfying all the usual properties, including a projection formula of the form
$$tr^G_H i^G_H (x)=[G:H]x,$$as well as the double coset formula. Thus, or course, $h^*G=H^*(BG,\F_2)$ is a possibility, as well as $h^*G=\cht{BG}$ from \cite{pedro2}, appendix. We could also take $h^*G=\textrm{Chern}^*G$, the smallest subring of $H^*(BG,\F_2)$ containing all Chern classes and transfers of such, or the analogous subring of $\cht{BG}$.

\begin{thm}[C\'ardenas-Kuhn] Suppose $L\subset K \subset G$ is a closed system, with $L$ elementary $2$-abelian. Suppose also that $W_K(L)$ contains a $2$-Sylow of $W_G(L)$. Then
$$im(h^*G \to h^*L)=(h^*L)^{W_G(L)} \cap im(h^*K \to h^*L).$$
\end{thm}

A proof can be found in \cite{adem} in the case of mod $2$ cohomology. It can be formally transposed to the other choices of $h^*$ as above.

\subsection{The mod $2$ Chow ring of $G_2(q)$.}

\begin{thm}\label{thm:g2q} Suppose $q \equiv 1$ mod $4$. Then the Chow ring $\cht{BG_2(q)}$ contains a polynomial subring $\F_2[d_4,d_6,d_7]$, which injects into the mod $2$ cohomology ring. Moreover, the classes $d_4$, $d_6$, $d_7$ are transfers of Chern classes.
\end{thm}

\begin{proof} The group $G_2(q)$ has two conjugacy classes of elementary abelian $2$-subgroups, represented, say, by $I$ and $II$ where the normalizer $E$ of $I$ fits in the exact sequence
$$\begin{CD} 1 @>>> I @>>> E @>>> GL_3(\F_2) @>>> 1 \end{CD}.$$
See \cite{m1}, 4.6, for all this. Moreover, $E$ contains a $2$-Sylow subgroup $S$ of $G_2(q)$, and it follows from \cite{m2}, 2.4 and 7.1, that $I\subset S \subset G_2(q)$ and $I\subset S \subset E$ are closed systems.

The invariants to consider are Dickson invariants:
$$(H^*(BI,\F_2))^{GL_3(\F_2)}=\F_2[D_4,D_6,D_7], \quad (\cht{BI})^{GL_3(\F_2)}=\F_2[d_4,d_6,d_7]$$ and $d_i$ maps to $D_i^2$ via the cycle map. Applying the C\'ardenas-Kuhn theorem twice, we are left with the task of proving that $\cht{BE}$ contains classes $d_4$, $d_6$, and $d_7$ mapping to the elements with the same name in $\cht{BI}$.

To this effect, we use the fact (\cite{m2},\S 9) that $E$ can be embedded as a subgroup of $G_2$, in such a way that the composition
$$H^*(BG_2,\F_2) \to H^*(BE,\F_2) \to H^*(BI,\F_2)$$ is injective and sends the Siefel-Whitney class $w_i$ as in the introduction to $D_i$. It follows that the Chern class $c_i\in \cht{B\gc}$ restricts to $d_i \in \cht{BI}$, and we are done.

To prove the last statement, one only has to run through the argument above replacing throughout the Chow ring by its smallest subring containing all Chern classes and transfers of such. Since the C\'ardenas-Kuhn theorem applies to this theory, and since the Chow ring of an elementary abelian group is generated by Chern classes, the result follows.
\end{proof}

\begin{coro} Unlike the situation when $p$ is odd, the restriction 
$$\cht{BG_2(q)} \to \cht{BT(q)}$$ is not injective.
\end{coro}

\begin{proof} Compare Krull dimensions (or transcendance degrees). \end{proof}

\pagebreak
\begin{center}
{\sc\Large Part II: The Chow ring of $\Spin_7\c$}
\end{center}

\section{Actions of $G_2$ and $Spin(7)$}\label{section:prelim}

\subsection{Complexifications.} We shall need a few precise results on complexifications (\cite{btd}, III, 8). If $G$ is a compact Lie groups, we denote its complexification by $G_\C$. Recall that any representation $G \to \GL_n\c$ extends to $G_\C$. More precisely, we have the following proposition, in which $\td{A}=~^t\bar{A}^{-1}$.

\begin{prop} Let $r: G \to \GL_n\c$ be a faithful, real representation. Then $r$ extends to a faithful representation of $G_\C$. Moreover, $G=U(n)\cap G_\C$.

Conversely, if $G_\C$ is any algebraic subgroup of a $\GL_n\c$, and if $G_\C$ is stable under the operation $A\mapsto \td{A}$, then $G_\C$ is a complexification of the compact Lie group $G=U(n)\cap G_\C$.
\end{prop}

See {\em loc. cit.} and exercise 8.8.7 there in particular.

\subsection{Actions of $\gc$.} The action of the compact group $G_2$ on $\VR$ preserves norms, and it will be very useful to us that $G_2$ acts transitively on pairs of orthogonal unit vectors: see \cite{adams}.

The vector space $\VC$ carries the standard quadratic form: $$q(z_1,\ldots,z_7)=z_1^2 + \cdots + z_7^2$$ which is obtained by complexifying the standard inner product of $\VR=\R^7$. This form is preserved by the action of $\gc$. We let $Q=\{\z \in \VC : q(\z)=1 \}$.

\begin{prop}\label{prop:qtrans} The action of $\gc$ on $Q$ is transitive. Moreover, the stabilizer of any point is an $\SL_3\c$. 
\end{prop}

\begin{proof} The action of $\O_7\c$ on $Q$ is transitive. It follows that $Q$ is a smooth, connected variety of (complex) dimension $6$. We prove that the action of $\gc$ is transitive by showing that all the orbits are open.

Consider an element $\z=\x + i\y \in \VC$ with $\x, \y \in \VR$. The condition $q(\z)=1$ amounts precisely to $\| \x \|=\| \y \| + 1$ (Euclidean norms) and $\x \perp \y$.

Suppose first that $\y\ne 0$. Since the compact group $G_2$ acts transitively on pairs of orthogonal vectors in $S^6$, it follows that the orbit of $\z$ contains a real submanifold diffeomorphic to the unit tangent bundle of $S^6$, so it has real dimension $\ge 11$. But the orbit of $\z$ is a complex manifold and must have even real dimension. It follows that this orbit has real dimension $12$ (complex dimension $6$), so it is open in $Q$.
 
Now suppose that $\y=0$, so that $\z\in\VR$. We prove that the stabilizer of $\z$ is isomorphic to $\SL_3\c$. Since $\dim \gc=14$ and $\dim \SL_3\c=8$, it will follow that the orbit of $\z$ has dimension $6$ and is therefore open in $Q$, so the proof will be complete. 
 
We may as well take $\z=e_1$, first vector in the standard basis of $\VR$. The stabilizer of $\z$ under the action of $O_7\c$ is $O_6\c$ embedded in the obvious way. Since $S=\O_6\c \cap \gc$ is stable under the operation $A\mapsto \td{A}$, it follows from the proposition above that $S$ is a complexification of $S\cap U(7)= O(6) \cap G_2$; in other words, the stabilizer in $\gc$ is a complexification of the stabilizer in $G_2$. But the latter is an $SU(3)$.
\end{proof}

Now we let $C=\{\z\ne 0 \in \VC : q(\z)=0 \}$.

\begin{prop}\label{prop:ctrans} The action of $\gc$ on $C$ is transitive, and the stabilizer of any point is a semi-direct product $H\rtimes \SL_2\c$, where $H$ is a $1$-connected solvable group.
\end{prop}

\begin{proof} Again, the action of $\O_7\c$ is transitive, so that $C$ is a smooth, connected variety. In fact, $\z = \x +i\y \in C$ if and only if $\|\x\|=\|\y\|$ and $\x\perp \y$, and it is easy to see that $C$ is diffeomorphic to $TS^6 \times \R^{>0}$ where $TS^6$ is the unit tangent bundle of $S^6$. In particular, $C$ is $4$-connected.

We may argue as in the preceding proof to see that the orbits under $\gc$ are open, so that the action is transitive. Let $S$ denote the stabilizer of $e_1$, the first vector in the standard basis of $\VR$. 

The group $S$ is simply-connected. Let $H$ denote its radical, the maximal normal, solvable, connected subgroup. Then $S/H$ is a semi-simple complex Lie group, and it is therefore the complexification of a compact Lie group $G$. It follows that $\pi_2(S/H)=0$, so that $H$ is simply-connected.

The group $S/H$ (or equivalently, the group $G$) is simply-connected, has rank less than  or equal to $2$ (which is the rank of $G_2$), and dimension less or equal to $8$ (since $C$ has complex dimension $6$). From the classification, $S/H$ is either $\SL_2\c$, $\SL_2\c \times \SL_2\c$, or $\SL_3\c$.

The case $\SL_3\c$ is impossible, since this group has dimension $8$: it would follow that $S=\SL_3\c$ and that $Q$ and $C$ are isomorphic algebraic varieties. However $Q=\O_7\c/\O_6\c$ has the homotopy groups of $S^6$ and these do not fit into the exact sequence of homotopy groups associated to the fibration $S^5 \to TS^6 \to S^6$. Thus $Q$ and $C$ are not even homotopy equivalent.

We shall rule out the case $\SL_2\c \times \SL_2\c$ as well. Since $H$ is solvable, it is made up of succesive extensions by either $\ga$ or $\gm$; and since $\pi_1(H)=0$, it is easy to see that $H$ is isomorphic to affine space as a variety. Hence $S$ has the homotopy type of $S/H$. Consider the Serre spectral sequence associated to the fibration $S\to \gc \to C$. This has $E_2^{p,q}=H^p(C, H^q(S,\Q))$ and converges to $H^{p+q}(\gc,\Q)$. As $C$ is $4$-connected, all differentials vanish on $E_r^{0,3}$ and we conclude that $H^3(S,\Q)=H^3(\gc,\Q)=\Q$. As a result, $S/H$ can only be an $\SL_2\c$.

It remains to see that the extension $H\to S \to S/H$ is split. Indeed, the action of $G_2$ on pairs of orthogonal vectors in $S^6$ having stabilizers $SU(2)$, the action of $\gc$ must contain a copy of the complexification $\SL_2\c$ in every stabilizer. This group is semi-simple and must have trivial intersection with the radical.\end{proof}

%

\subsection{Actions of $Spin(7)$.} We turn to $Spin(7)$ and its action on the spin representation $\Delta_\R$. We write $\DC$ or simply $\Delta$ for $\Delta_\R\otimes_\R \C$. Let us fix a basis $e_1,\ldots,e_8$ of $\DR$ such that the first seven of these vectors span a subspace isomorphic to $\VR$ as a representation of $G_2\subset Spin(7)$ (while $e_8$ is fixed by $G_2$). As before, $\Spin_7\c$ preserves a quadratic form on $\DC$, and we can arrange for it to be the standard one by choosing the basis appropriately.

We still write $q$ for this form. We define $Q^7$ and $C^7$ by analogy with $Q$ and $C$ above:
$$Q^7=\{\z \in \DC : q(\z)=1 \}, \quad C^7=\{\z\ne 0 \in \DC : q(\z)=0 \}.$$
These are algebraic varieties in the $8$-dimensional affine space $\DC$, and the upperscript gives the complex dimension. The propositions which follow are completely analogous to those in the previous section, and the proofs will be a little sketchy.

\begin{prop}\label{prop:qtransspin} The action of $\Spin_7\c$ on $Q^7$ is transitive. Moreover, the stabilizer of any point is isomorphic to $\gc$. 
\end{prop}

\begin{prop}\label{prop:ctransspin} The action of $\Spin_7\c$ on $C^7$ is transitive, and the stabilizer of any point is a semi-direct product $H\rtimes \SL_3\c$, where $H$ is a $1$-connected solvable group.
\end{prop}

\begin{proof}[Proof of \ref{prop:qtransspin}] We proceed as in \ref{prop:qtrans}. Again, we show that the action is transitive by showing that the orbits are open in the smooth, connected variety $Q^7$.

Pick $\z\in Q^7$, $\z=\x + i\y \in \DC$ with $\x, \y \in \DR$. Since the action of $Spin(7)$ on $S^7\subset\DR$ is transitive, we may assume that either $\y\ne 0$ or $\z=e_8$. In the former case we use the fact that the action of $Spin(7)$ is transitive on pairs of orthogonal unit vectors in $\DR$; in the latter case we argue that the stabilizer of $e_8$ has to be the complexification of $G_2$. In either case this is completely analogous to \ref{prop:qtrans}.
%
\end{proof}

\begin{proof}[Proof of \ref{prop:ctransspin}] We leave it to the reader to prove that the action is transitive.

Let $S$ be the stabilizer of a point $\z$ with $q(\z)=0$. By proceeding exactly as in \ref{prop:ctrans}, things come down to identifying a semi-simple, $1$-connected algebraic group $S/H$ as $\SL_3\c$. In fact, $S$ is simply-connected, and as a result  its radical $H$ is simply-connected as well, so it does not contain a torus, and we may find a Levi subgroup $L$ such that $S=H\rtimes L$.

Let us prove first that $S$, and thus $S/H$ or $L$, contains a copy of $\SL_3\c$. Indeed, consider $\z=(1+i)e_1 + (1-i)e_8$: this is stabilized at least by the copy of $SU(3)$ in the compact group $G_2$ which fixes $e_1$, and the claim follows.

Thus $L$ has rank between $2$ and $3$ and dimension between $8$ and $14$ (since the dimension of $\Spin_7\c$ is $21$ and $\Spin_7\c/S$ has dimension $7$). Suppose that $L$ has rank $3$. The reductive subgroups of maximal rank in $\Spin_7\c$ are known up to local isomorphism (\cite{borel2}): either $L$ has the type of $\Spin(5)\times \gm$ (which is not semi-simple), or that of $\Spin(6)$ (which has dimension $15$), or finally it could have the type of $(\SL_2\c)^3$ (which is not simple and thus gives the wrong rational cohomology, as in \ref{prop:ctrans}). Therefore $L$ must have rank $2$ and contain $\SL_3\c$.

We conclude that $L$ is either $\SL_3\c$ or $\gc$; but $L=\gc=S$ would imply that $Q^7$ and $C^7$ are diffeomorphic, which they are not.
\end{proof}

\section{Stratification of $\Delta$.}\label{section:stratspin}

In this section, we implement the "stratification method", which was first used by Vezzosi \cite{vezzosi}, and later by Vistoli \cite{vistoli2} and Molina-Vistoli \cite{vistoli1}. 

The idea is essentially as follows. One makes use of the equivariant Chow rings $CH^*_GX$ introduced by Edidin-Graham in \cite{didingra} buiding on Totaro's construction in \cite{totaro}, for any algebraic group $G$ and any variety $X$ endowed with an action of $G$; when $X$ is reduced to a point, $CH^*_G(point)=CH^*BG$. These have the same formal properties as ordinary Chow groups. 

The starting point of the stratification method is the fact that $CH^*_G W=CH^*BG$ for any representation $W$ of $G$ (which is a $G$-equivariant vector bundle over a point, hence the isomorphism). Then, one uses the localisation exact sequence:
$$\begin{CD} CH^{*-c}_G Y @>>> CH^*_G W @>>> CH^*_G (W-Y) @>>> 0\end{CD}$$ which one gets for any $G$-invariant closed subvariety $Y$ of codimension $c$ in $W$. For instance, if the action is transitive on both $Y=G/H_1$ and $W-Y=G/H_2$, then we can gain precious information on $CH^*BG$ by using the elementary isomorphism $CH^*_G(G/H_i)=CH^*BH_i$: hopefully the subgroups $H_i$, which are smaller than $G$, may each have a Chow ring which is already known.

However, the action will typically not be transitive on $W-Y$, so one has to look for a closed orbit and repeat the process: hence the term "stratification".

We shall now carry this out in the case of $\Spin_7\c$ and its spin representation $\D$. We follow rather closely the method used in $\cite{vistoli1}$ for the groups $\O_n$.

\medskip
\noindent{\em Notations.} In what follows, we shall write $c_i(\D)$, resp. $c_i(V)$, for the $i$-th Chern class of $\DC$, resp. $\VC$, in the Chow ring of $\Spin_7\c$. They both restrict to the element $c_i$ in the Chow ring of $\gc$.
\medskip

\subsection{Two exact sequences.}

The first arrow in the localization sequence:
$$\begin{CD}CH^{*-7}_{\Spin_7\c} \{0\} @>>> CH^*_{\Spin_7\c} \DC @>>> CH^*_{\Spin_7\c}(\DC-\{0\}) @>>> 0 \end{CD}$$ is multiplication by $c_8(\D)$, so we draw
$$CH^*_{\Spin_7\c}(\DC-\{0\})=\frac{CH^*B\Spin_7\c}{(c_8(\D))}.$$
This fits into another localization sequence, where $U=\DC - (C^7\cup \{0\})$: 
$$\begin{CD}CH^{*-1}_{\Spin_7\c} C^7 @>>> CH^*_{\Spin_7\c} (\DC -\{0\}) @>>> CH^*_{\Spin_7\c} U @>>> 0 \end{CD}.$$
 Our task is to identify the terms in this sequence.

\subsection{The term $CH^*_{\Spin_7\c} U$.} Consider first the map $p: \gm \times Q^7 \to U$ given by $(t,\z)\mapsto t\z$. There is an action of $\Z/2 \times \Spin_7\c$ on $\gm\times Q^7$ obtained by extending the action of $\Spin_7\c$ trivially on $\gm$ and setting $\epsilon\cdot (t,\z)=(- t, -\z)$, where $\epsilon$ generates the factor $\Z/2$. The action of $\Z/2$ alone is free and the quotient is precisely $U$, with quotient map $p$; therefore we have 
\begin{equation} CH^*_{\Spin_7\c} U=CH^*_{\Z/2 \times \Spin_7\c} (\gm \times Q^7). \end{equation}

We may extend the action of $\Z/2\times \gc$ to $\af^1\times Q^7$ in the obvious way and then consider the exact sequence:

$$\begin{CD}CH^{*-1}_{\Z/2\times \Spin_7\c} (\{0\}\times Q^7) @>>> CH^*_{\Z/2\times \Spin_7\c}(\af^1\times Q^7)  @>>> CH^*_{\Z/2\times \Spin_7\c} (\gm\times Q^7) @>>> 0 \end{CD}$$

Again the first arrow here is multiplication by the top Chern class of the representation obtained by composing the projection $\Z/2 \times \Spin_7\c \to \Z/2$ and the non-trivial character of $\Z/2$. Calling it $\alpha$, it follows that 
\begin{equation} CH^*_{\Spin_7\c} U=\frac{CH^*_{\Z/2 \times \Spin_7\c} Q^7}{(\alpha)}.\end{equation}

It follows from proposition \ref{prop:qtransspin} that the action of $\Z/2\times \Spin_7\c$ on $Q^7$ is transitive. The stabilizers are isomorphic to $\Z/2 \times \gc$ embedded in $\Z/2 \times \Spin_7\c$ thus: $$(1,g)\mapsto (1,g), \quad (\epsilon,g)\mapsto (\epsilon,\tau g)$$ where $\tau$ generates the kernel of $\Spin_7\c\to \SO_7\c$. We note:
$$\begin{array}{rcl} 
CH^*_{\Z/2 \times \Spin_7\c} Q^7 & = & CH^*_{\Z/2 \times \Spin_7\c} (\Z/2 \times \Spin_7\c)/(\Z/2 \times \gc)
\\                              & = & CH^* B(\Z/2\times \gc)
\\                              & = & CH^*B\gc[\alpha]. 
\end{array} $$ so that, finally, $CH^*_{\Spin_7\c} U= CH^*B\gc$.

\subsection{The term $CH^*_{\Spin_7\c} C^7$.} From proposition \ref{prop:ctransspin}, the action of $\Spin_7\c$ on $C^7$ is transitive with stabilizer $H\rtimes \SL_3\c$. The group $H$ is solvable and simply connected, so as a variety it is an affine space. Therefore the map $B\SL_3\c \to B(H\rtimes \SL_3\c)$ induced by inclusion has fibres isomorphic to affine space, and thus $CH^*B(H\rtimes \SL_3\c)=CH^*B\SL_3\c$. So we draw
$$CH^*_{\Spin_7\c} C^7=CH^*B\SL_3\c=\Z[x_2,x_3].$$

\subsection{Conclusion.}\label{subsec:conclustrat} Thus far, the exact sequence in the introduction to this section has been identified as:
$$\begin{CD}CH^{*-1} B\SL_3\c @>{i_*}>> CH^* B\Spin_7\c /(c_8(\D)) @>{j^*}>> CH^* B\gc @>>> 0 \end{CD}.$$

\begin{rmk}\label{rmk:projmodc8} A couple of comments are in order here. First, we shall let the reader unwind all the definitions and verify that $j^*$ may be identified with the restriction map induced by the inclusion $\gc \to \Spin_7\c$. Likewise, $i^*$ can be seen as the map induced by the inclusion $\SL_3\c \to \Spin_7\c$ (notice that $c_8(\D)$ restricts to $0$ in $CH^*B\SL_3\c$, see below). However, it would not be proper to call $i_*$ the "transfer" from $\SL_3\c$: indeed, this subgroup does not have finite index in $\Spin_7\c$.

We nevertheless have a projection formula between $i^*$ and $i_*$. We shall use it in $CH^*B\Spin_7\c$ in the form $i_*(xi^*(y))=i_*(x)y \mod (c_8(\D)).$
\end{rmk}

Now, the representation $\DC$ splits up, as an $\SL_3\c$-module, as the direct sum of a trivial module and $W\oplus W^*$, where $W$ is the standard representation of $\SL_3\c$. From this we draw $i^*(c_2)=2x_2$, $i^*(c_4)=x_2^2$, $i^*(c_6)=x_3^2$, and all the other Chern classes restrict to $0$. The ring $CH^*B\SL_3\c=\Z[x_2,x_3]$ is generated as a module over $\Z[x_2^2,x_3^2]$ by $1$, $x_2$, $x_3$ and $x_2x_3$. If we put $\zeta_1=i_*(1)$, $\zeta_3=i_*(x_2)$, $\zeta_4=i_*(x_3)$ and $\zeta_6=i_*(x_2x_3)$, we conclude from the projection formula that the image of $i_*$ is the ideal generated by $\zeta_1$, $\zeta_3$, $\zeta_4$, $\zeta_6$. (Note that these elements have degree less than that of $c_8(\D)$, so we may safely regard them as living in $CH^*B\Spin_7\c$.) Moreover we have:

\begin{lem}\label{lem:zeta1} $\zeta_1=0$.
\end{lem}

\begin{proof}(See \cite{vistoli1}.) The element $i_*(1)$ corresponds to $[C^7]\in CH^1_{\Spin_7\c} (\D - \{ 0 \})$, and $C^7=q^{-1}(0)$ by definition. However, the equivariant map $q:\D - \{ 0 \} \to \af^1$ is flat, so $[q^{-1}(0)]=q^*([0])$. Finally, note that $[0]\in CH^1_{\Spin_7\c} \af^1$ is the top Chern class of the trivial representation, so it is $0$.
\end{proof}

In particular, it follows from this and the relation $i^*(c_2)=2x_2$ that $2\zeta_3=0$. Let us summarize the information gathered in this section:

\begin{prop}\label{prop:sping2} Let $c_i$ denote either $c_i(\D)$ or $c_i(V)$, for $1\le i \le 7$. The Chow ring of $\Spin_7\c$ is generated by the classes $c_2$, $c_4$, $c_6$, $c_7$, $c_8(\D)$ and $\zeta_3$, $\zeta_4$, $\zeta_6$. We have the relation $2\zeta_3=0$ and $$CH^*B\Spin_7\c/(\zeta_3, \zeta_4, \zeta_6, c_8(\D))=CH^*B\gc.$$
\end{prop}

\begin{rmk}\label{rmk:kerasmod} Consider the subring $B$ generated by $c_4(V)$ and $c_6(V)$. From the computations above, we notice that the kernel of the restriction map from $CH^*B\Spin_7\c/(c_8(\D))$ to $CH^*B\gc$ is the submodule generated over $B$ by $\zeta_3$, $\zeta_4$, and $\zeta_6$.
\end{rmk}

\section{Input from Brown-Peterson cohomology.}\label{section:otherinput}

\subsection{Torsion.} 
Let $\alpha\in Spin(7)$ be the element of order two giving "the outer
automorphism of $SO(6)$ which becomes inner in $SO(7)$" as in
\cite{adams}, p34. (In standard notations on Clifford algebras,
$\alpha=e_6e_7$.) This element normalizes the copy of $\Spin_6\c$
contained in $\Spin_7\c$, and together they generate a subgroup
$\Spin_6\c \rtimes \Z/2$ which contains the full normaliser of a
maximal torus in $\Spin_7\c$. Since $\Spin_6\c=\SL_4\c$ has a
torsion-free Chow ring, we conclude as before using theorem
\ref{thm:totaro} that

\begin{lem}\label{lem:torsionspin} If a class $x\in CH^*B\Spin_7\c$ is torsion, then $2x=0$.
\end{lem}

Therefore, there is little loss of information in considering
$$CH^*(B\Spin_7\c)_{(2)}=CH^*B\Spin_7\c \otimes_\Z \Z_{(2)}$$ instead
of the integral Chow ring (indeed, additively, there is none).

\begin{rmk} 
The lemma above also follows from general principles. That is, torsion
elements in the Chow ring of an algebraic group $G$ are always killed
by the "torsion index", which in the case of $\Spin_7\c$ is $2$: see
\cite{totarospin} for details.
\end{rmk}

\subsection{Brown-Peterson cohomology.} 
We shall use the natural map defined for any complex algebraic variety
by Totaro (see \cite{totaro} and the references there):
$$CH^*X \to BP^*X \otimes_{BP^*} \Z_{(2)}$$ where $BP$ denotes the
Brown-Peterson spectrum at the prime $2$. Recall that
$$BP^*(point)=\Z_{(2)}[v_1,v_2,\ldots]$$ with $v_i$ of degree $-2(2^i
-1)$.

The $BP$-cohomology of $B\Spin_7\c$ has been studied by Kono and
Yagita, see \cite{konoyagita} and also \cite{yagshu} . To state their
result, we need to introduce some notations. If $R$ is a ring and the
$x_i$'s are indeterminates, we write $R[x_i]$ for the ring of
polynomials with coefficients in $R$, while $R\<x_i\>$ denotes the
free $R$-module with basis given by the $x_i$'s; if the $x_i$'s live
in an $R$-algebra $S$, resp. an $R$-module $M$, the notations imply
that they are algebraically, resp. linearly, independant over $R$. We
use the notation $R\{x_i\}$ for the (non necessarily free) submodule
of $M$ generated over $R$ by the $x_i$'s. Thus for example,
$R[x]\{x\}=R[x]\<x\>$ means the (non unitary) ring of polynomials
without constant terms.

We let $E_r^{*,*}$ be the $r$-th page of the Atiyah-Hirzebruch
spectral sequence
$$H^*(B\Spin_7\c, \Z_{(2)})\otimes_{\Z_{(2)}} BP^* \Rightarrow
BP^*B\Spin_7\c.$$ Finally, the letter $w_i$ will denote any element in
$H^*(B\Spin_7\c, \Z_{(2)})$ whose reduction mod $2$ is the
Siefel-Whitney class which we also denote by $w_i$ (we take it as part
of the result that follows that $w_i$ can be lifted this way). Note
that $w_8$ is the top Stiefel-Whitney class of $\D$, while the others
refer to $\VR$.

Then Kono and Yagita prove the following:
$$E_\infty^{*,*}=A\{1,2w_4,v_1w_8,2w_8,2w_4w_8\} \bigoplus A[w_7^2]\{
w_7^2 \}/(2,v_1,v_2,v_3w_7^2w_8^2)$$ where
$A=BP^*[w_4^2,w_6^2,w_8^2]$. Here subscripts give the cohomological
degree. Note that the first summand is isomorphic to
$$A\<x_0,x_4,x_6,x_8,x_{12}\>/(2x_6-v_1x_8).$$
As a result, we have the following additive isomorphism:
$$BP^*B\Spin_7\c \otimes_{BP^*} \Z_{(2)}=A'\<\tilde x_0,\tilde
x_4,\tilde x_6,\tilde x_8,\tilde x_{12}\>/(2\tilde x_6) \bigoplus
A'/(2)[\tilde w_7^2]\{ \tilde w_7^2 \} $$ where $A'=\Z_{(2)}[\tilde
w_4^2,\tilde w_6^2,\tilde w_8^2]$, and where $\tilde x_0$, $\tilde
x_4$, etc, is any class in $BP^*B\Spin_7\c$ which is represented on
the page $E_\infty^{*,*}$ by $x_0=1$, $x_4=2w_4$, etc. To prove the
relation $2\tilde x_6=0$, we note that $BP^*B\Spin_7\c \otimes_{BP^*}
\Q= H^*(B\Spin_7\c,\Q)$, and this is a polynomial ring in the
Pontryagin classes $c_2(V)$, $c_4(V)$, $c_6(V)$, so $\tilde x_6$ must
be torsion. However, lemma \ref{lem:torsionspin} applies here, since
it relies solely on theorem \ref{thm:totaro}, and in turn this result
has an obvious analog for any cohomology $h^*$: indeed $BG$ is a
stable retract of $BN$, so $h^*(BG)$ is a direct summand in $h^*(BN)$.

There are several ways of picking the $\tilde x_i$'s, and some choices
are more natural than others: for example if we have a class in
$BP^*B\Spin_7\c$ which under the cycle map to
$H^*(B\Spin_7\c,\Z_{(2)})$ gives $2w_4$, we may take it for $\tilde
x_4$. This is clearly not possible for $v_1w_8$, but it is feasible in
all the other cases: indeed Yagita has shown that we could take
$c_2(\D)$, $c_4(\D)$, $c_6(\D)$ for $\tilde x_4$, $\tilde x_8$,
$\tilde x_{12}$ respectively. Moreover, we can take $c_4(V)$,
$c_6(V)$, $c_7(V)$ and $c_8(\D)$ for $\tilde w_4^2$, $\tilde w_6^2$,
$\tilde w_7^2$ and $\tilde w_8^2$.

Thus we obtain the following proposition. Before stating it, we introduce some 

\smallskip
\noindent{\em Notations.} From now on, whenever it will make a formula
easier to read, we shall write $c_i$ for $c_i(V)$ and $c'_i$ for
$c_i(\D)$.

\begin{prop}[Yagita]\label{prop:yagita} The cycle map induces an epimorphism:
$$CH^*(B\Spin_7\c)_{(2)} \to \Z_{(2)}[c_4,c_6,c'_8]\otimes \Big( \Z_{(2)}\<1,c'_2,c'_4,c'_6\> \oplus \F_2\<\tilde x_6\> \oplus \F_2[c_7]\<c_7\> \Big).$$
\end{prop}

\begin{coro}\label{coro:generators} We have $$c_4(\D) - c_4(V) = a\zeta_4, \qquad c_6(\D) - c_6(V) = b\zeta_6$$ where $a$ and $b$ are odd integers. As a result, the ring $CH^*(B\Spin_7\c)_{(2)}$ is generated by the classes $c_2(\D)$, $c_4(\D)$, $c_4(V)$, $c_6(\D)$, $c_6(V)$, $c_7(V)$, $c_8(\D)$ and $\zeta_3$. The kernel of the restriction map to the Chow ring of $\gc$ is the ideal generated by $c_8(\D)$, $\zeta_3$, $c_4(\D) - c_4(V)$, and $c_6(\D) - c_6(V)$. In the notations of the Proposition, $\zeta_3$ maps to $\tilde x_6$.
\end{coro}

\begin{proof} The class $c_4(\D)-c_4(V)$ restricts to $0$ in the Chow ring of $\gc$, and we deduce from the results of the previous section that $c_4(\D) - c_4(V) = a\zeta_4$ for an integer $a$. The last proposition implies then that $a$ is odd. The argument for $c_6(\D) - c_6(V)$ is similar, and the rest is clear.
\end{proof}

\begin{rmk2} In the course of the following proof, we shall show that $a=3$ and $b=1$.
\end{rmk2}

\section{Final steps}\label{section:end}

\subsection{Multiplicative relations.}

\begin{prop}\label{prop:relations} The following relations hold in $CH^*(B\Spin_7\c)_{(2)}$, where $\delta_i=0$ or $1$ ($i=1,2)$.
$$\begin{array}{rlrl}
(1) & \zeta_3^2=0 & (8) & c'_2c'_6 - c'_2c_6=\frac{2}{3}c_4(c'_4 - c_4) + 16c'_8  
\\ (2) & \zeta_3c_7=0 &  (9) & c'_2c_7=\delta_1c_6\zeta_3
\\ (3) & \zeta_3c'_4=\zeta_3c_4 & (10) &  c'_4(c'_4 - c_4)=c_4(c'_4 - c_4) + 36 c'_8  
\\ (4) & \zeta_3c'_6=\zeta_3c_6 & (11) & c'_4(c'_6 - c_6)=c_4(c'_6 - c_6) +6c'_2c'_8  
\\ (5) & \zeta_3c'_2=0 & (12) & c'_4c_7 - c_4c_7=\delta_2c_8\zeta_3 
\\ (6) & (c'_2)^2 - 4c_4=\frac{8}{3}(c_4'-c_4) & (13) & c'_6(c'_6 - c_6)=c_6(c'_6 - c_6) + c'_8(\frac{8}{3}c'_4 + \frac{4}{3}c_4)
\\ (7) & c'_2c'_4 - c'_2c_4 = 6(c'_6 - c_6) & (14) & c'_6c_7 - c_6c_7=0
\end{array}$$
\end{prop}

\begin{proof} We use the notations of \S\ref{section:stratspin}:
$$\begin{CD} CH^*B\Spin_7\c @>{i^*}>> CH^*B\SL_3\c=\Z[x_2, x_3]\end{CD}.$$
Consider, say, equation (3). We have
$$\zeta_3c'_4=i_*(x_2)c'_4=i_*(x_2i^*(c'_4))=i_*(x_2i^*(c_4))=i_*(x_2)c_4=\zeta_3c_4.$$
If we proceed in the exact same way, using the relations $i^*(c_2)=2x_2$, $i^*(c_4)=x_2^2$, $i^*(c_6)=x_3^2$, and $i^*(c_7)=0$, as well as the corresponding equations with $c'_i$ replacing $c_i$, and the last corollary, we see that equations (2), (4), (5), and (10) follow.

As for (1), it is only a matter of proving that $i^*(\zeta_3)=0$. However we certainly have $i^*(\zeta_3)=kx_3$ for an integer $k$, which is then a torsion number for $\zeta_4=i_*(x_3)$ (since $i_*(1)=0$, see \ref{lem:zeta1}). Consequently, $k=0$ from corollary \ref{coro:generators}, as $c'_4 - c_4$ is not torsion (cf \ref{prop:yagita}).

To continue, we need to use the fact that the projection formula is only valid modulo $c'_8$ (remark \ref{rmk:projmodc8}). From equation (5), we see that there are no elements of degree $5$ at all, so we have (14) at once. Otherwise we obtain:
$$\begin{array}{rlrl}
(7') & c'_2(c'_4 - c_4)=2a(c'_6 - c_6) & (11') & c'_4(c'_6 - c_6)=c_4(c'_6 - c_6) + \lambda_1c'_2c'_8
\\(8') & c'_2(c'_6 - c_6)=\frac{2b}{a} c_4(c'_4 - c_4) + \lambda_2c'_8 & (12) & c'_4c_7 - c_4c_7=\delta_2c_8\zeta_3 
\\(10') & c'_2(c'_4 - c_4)=c_4(c'_4 - c_4) + \lambda_3c'_8 & (13') & c'_6(c'_6 - c_6)=c_6(c'_6 - c_6) + c'_8(Ac'_4 + Bc_4)
\end{array}$$
To compute all the constants (except $\delta_2$), we investigate the restriction of the classes above to the Chow ring of a maximal torus. If we write $CH^*(BT)_{(2)}=\Z_{(2)}[x,y,z]$, we have
$$\sum_{i=0}^8 c'_iX^i=\prod_{\epsilon_i=\pm 1} (1 + (\epsilon_1 x + \epsilon_2 y + \epsilon_3 z)X)$$ while
$$\sum_{i=0}^7 c_iX^i=(1+2xX)(1-2xX)(1+2yX)(1-2yX)(1+2zX)(1-2zX).$$
Let us illustrate the method with equation (13), which is the hardest. We compute both sides of (13') modulo the ideal $(y,z)$ and we get $16x^{12}=6Ax^{12}$, so $A=\frac{8}{3}$. We replace $A$ by this value and compute equation (13') modulo $(z)$: we find that it holds if and only if we have
$$128 x^{6} y^{6} - {  \frac {256}{3}}  x^{
8} y^{4} + {  \frac {64}{3}}  x^{10} y^{2} + 
{  \frac {64}{3}}  y^{10} x^{2} - {  
\frac {256}{3}}  y^{8} x^{4}$$ 
$$ - 16 y^{10} B x^{2} - 16 x^{10
} B y^{2} 
  + 64 x^{8} y^{4} B + 64 y^{8} x^{4} B - 96 x
^{6} y^{6} B=0.$$
This forces $B=\frac{4}{3}$. The other cases are similar.

There remains to establish equations (6) and (9). Both follow from the fact that the left hand side restrict to $0$ in the Chow ring of $\gc$, and using remark \ref{rmk:kerasmod}. There is a constant in equation (6), which we compute as above.
\end{proof}

\begin{coro}\label{coro:addspin} There is an isomorphism $$CH^*(B\Spin_7\c)_{(2)}= \Z_{(2)}[c_4,c_6,c'_8]\otimes \Big( \Z_{(2)}\<1,c'_2,c'_4,c'_6\> \oplus \F_2\< \zeta_3 \> \oplus \F_2[c_7]\<c_7\> \Big).$$ As a consequence, all the products in $CH^*(B\Spin_7\c)_{(2)}$ are determined by \ref{prop:relations}. Moreover, there is an isomorphism
$$CH^*(B\Spin_7\c)_{(2)}=BP^*(B\Spin_7\c)\otimes_{BP^*} \Z_{(2)}.$$
\end{coro}

\begin{proof} The classes $c_4, c_6, c'_8$ generate a polynomial subring $R$ in $CH^*(B\Spin_7\c)_{(2)}$, from \ref{prop:yagita}.

Consider the classes $c'_2, c'_4, c'_6, \zeta_3$ and $c_7^k$, for $k\ge 1$. Proposition \ref{prop:relations} computes all pairwise products, as well as the squares of the first four of these. It follows from this and \ref{coro:generators} that they generate $CH^*(B\Spin_7\c)_{(2)}$ as an $R$-module.

From \ref{prop:yagita} again, we conclude that $c'_2, c'_4$ and $c'_6$ must generate a free $R$-module, while $\zeta_3$ and the powers of $c_7$ generate a free $R/(2)$-module.
\end{proof}

\begin{rmk} 
The Chow ring of $\Spin_7$ does not inject into the cohomology: for
instance, the element $\zeta_3$ maps to $0$. Moreover, this element is
not a transfer of a combination of Chern classes. These two facts,
which are intimately related, follow from \cite{yagshu} where the case
of $BP$ is considered. However, we point out that $\zeta_3$ is some
sort of generalised transfer of a Chern class from $\SL_3$, in the
sense of remark \ref{rmk:projmodc8}.
\end{rmk}

\begin{rmk} 
It is possible to apply the stratification method to $\Spin_{2n}$ and
its "natural" representation via $\SO_{2n}$ (rather than the spin
representation, as we have done in the present paper). The details are
strictly parallel to those in \cite{vistoli1} where the case of
$\SO_{2n}$ is treated. One obtains:
$$\begin{CD} CH^*B\Spin_{2n - 2} @>>> CH^*B\Spin_{2n} @>>>
    CH^*B\Spin_{2n - 1} @>>> 0.
\end{CD}$$
However, in practice it seems difficult to proceed further. Consider
for example $\Spin_8$. We know the Chow rings of $\Spin_7$ and
$\Spin_6=\SL_4$. The exact sequence above allows one to find
generators for $CH^*B\Spin_{2n}$ as a module over the ring of Chern
classes. Potentially one may draw sharper information from the
determination of the Brown-Peterson cohomology of $\Spin_8$ by Kono
and Yagita in \cite{konoyagita}. At this point though, getting one's
hands on the multiplicative relations is still a challenge. The
missing ingredient is corollary \ref{coro:generators}, which allows us
in the case of $\Spin_7$ to replace the obscure classes $\zeta_4$ and
$\zeta_6$ by Chern classes, thus getting some control over the
restrictions that these may have to various subgroups, and the
like. Another dissatisfying aspect of this computation is that it has
none of the symmetry which triality ought to impose.\footnote{added
  21/02/2008: Luis Alberto Molina Rojas in his PhD thesis has
  completed the computation for $\Spin_8$.}
\end{rmk}

\section{Appendix}

This appendix is not part of the published version of this article. It
appears separately as an addendum.

\subsection{The missing proof}

Proposition 8.3 is given a proof in \cite{yag2} (the
explicit reference was missing, for which I apologize). However, in
this paper Yagita relies on a result of Totaro's which has not been
published: namely, the claim that
\[ CH^*(B\Spin_7)_{(2)} \to \Z_{(2)}[c_4,c_6,c'_8]\otimes \Big(
\Z_{(2)}\<1,c'_2,c'_4,c'_6\> \oplus \F_2\<\tilde x_6\> \oplus
\F_2[c_7]\<c_7\> \Big) \]
is surjective depends on Totaro's argument to the effect that $\tilde
x_6$ is indeed in the image of the cycle map. We prove this
now. Recall that the target of the above map is $BP^* B\Spin_7 \otimes_{BP^*} \Z_{(2)}$.

\begin{lem}
There is no $2$-torsion in (cohomological) degree 6 in 
$$ BP^*B\gc \otimes_{BP^*} \Z_{(2)}.   $$
\end{lem}

This follows, for example, from the computations in
\cite{konoyagita}.

The idea now is simply to use the stratification method, as in \S
7. The reader will easily obtain the following exact sequence, where
$*$ is even:
$$ BP^{*-1} B\SL_3 \to BP^* B\Spin_7 / (c_8(\Delta)) 
\to BP^* B\gc \to BP^{odd} B\SL_3. $$
(The last term is a typical symptom of the stratification method used
with a generalized cohomology rather than Chow rings, see
\cite{vistoli1}.)

It is well-known that $BP^{odd} B\SL_3 = 0$, so one ends up with the
following commutative diagram:
$$\mspace{-125mu}\begin{CD}
CH^{*-1} B\SL_3 @>>> CH^* B\Spin_7 / (c_8(\Delta)) @>>> CH^* B\gc @>>>
0 \\
     @VVV                @VVV                              @VVV   \\
 BP^{*-1} B\SL_3\otimes_{BP^*} \Z_{(2)} @>>> BP^*
 B\Spin_7\otimes_{BP^*} \Z_{(2)} / (c_8(\Delta)) @>>>  BP^*
 B\gc\otimes_{BP^*} \Z_{(2)} @>>> 0
\end{CD}$$
The element $\tilde x_6$ is $2$-torsion, so the last lemma shows
that this element comes from $\SL_3$. The cycle map is an isomorphism
for $\SL_3$, so the diagram shows that $\tilde x_6$ is the image of
$\zeta_3\in CH^3 B\Spin_7$, QED.

\subsection{Totaro's conjecture}

We shall now expand on remark 9.3.

In \cite{totaro}, Totaro conjectures the following: {\em the Chow ring
  $CH^*BG$ of any finite group $G$ is generated by Chern classes and
  transfers of Chern classes}. He also establishes that $CH^*BG$ is
  generated by Chern classes for $*\le 2$.

However, in \cite{yagshu}, Schuster and Yagita have proved that the
element $\tilde x_6$ in the $BP$-cohomology of $\Spin_7$ as above is
not a transfer of Chern classes. It now follows that $\zeta_3$ is not
a transfer of Chern classes, either. 

What is more, Schuster and Yagita also prove that the restriction of
$\tilde x_6$ to a certain finite subgroup $G$ of $\Spin_7$ is an
element in the $BP$-cohomology of $BG$ which, again, is not a transfer
of Chern classes. We see now that this element comes from $CH^3 BG$,
and gives an optimal counterexample to Totaro's conjecture.

\bibliography{myrefs}
\bibliographystyle{siam}
\end{document}